\documentclass[a4paper, 11pt] {article}
\usepackage[cp1251]{inputenc}
\usepackage[russian, english]{babel}
\usepackage{amsfonts}
\usepackage{mathtext}
\usepackage{cite}

\usepackage{graphicx}

\usepackage{verbatim}
\usepackage{amsmath}
\usepackage{amsthm}
\usepackage{amssymb}
\usepackage{delarray}
\usepackage{mathrsfs}
\usepackage{xcolor}

\begin{document}
%\vspace*{5mm}

%\thispagestyle{empty}

\begin{center}
{\bf\large ISO-BISPECTRAL POTENTIALS FOR STURM--LIOUVILLE-TYPE OPERATORS WITH SMALL DELAY}
\end{center}

\begin{center}
{\bf\large Neboj\v{s}a Djuri\'c\footnote{Faculty of Architecture, Civil Engineering and Geodesy, University
of Banja Luka, {\it nebojsa.djuric@aggf.unibl.org}} and Sergey Buterin\footnote{Department of Mathematics,
Saratov State University, {\it buterinsa@info.sgu.ru}}}
\end{center}

{\bf Abstract.} In recent years, there appeared a considerable interest in the inverse spectral theory for
functional-differential operators with constant delay. In particular, it is well known that, for each fixed
$\nu\in\{0,1\},$ the spectra of two operators generated by one and the expression $-y''(x)+q(x)y(x-a)$ and
the boundary conditions $y^{(\nu)}(0)=y^{(j)}(\pi)=0,$ $j=0,1,$ uniquely determine the complex-valued
square-integrable potential $q(x)$ vanishing on $(0,a)$ as soon as $a\in[\pi/2,\pi).$ For $a<\pi/2,$ the main
equation of the corresponding inverse problem is nonlinear, and it actually became the {\it basic question}
of the inverse spectral theory for Sturm--Liouville operators with constant delay whether the uniqueness
holds also in this nonlinear case. A few years ago, a positive answer was obtained for $a\in[2\pi/5,\pi/2).$
Recently, the authors gave, however, a {\it negative} answer for $a\in[\pi/3,2\pi/5)$ by constructing
infinite families of iso-bispectral potentials. Meanwhile, the question remained open for the most difficult
nonlinear case $a\in(0,\pi/3),$ allowing the parameter $a$ to approach the classical situation $a=0,$ in
which the uniqueness is well known. In the present paper, we address this gap and give a negative answer in
this remarkable case by constructing appropriate iso-bispectral potentials.

%\smallskip
{\it Key words:} Sturm--Liouville operator with small delay, functional-differential operator, inverse
spectral problem, iso-bispectral potentials

%\smallskip
{\it 2010 Mathematics Subject Classification:} 34A55 34K29
\\

{\large\bf 1. Introduction}
\\

Inverse problems of spectral analysis consist in recovering operators from their spectral characteristics.
One of the first results in the inverse spectral theory says that the spectra of two boundary value problems
for one and the same Sturm--Liouville equation with one common boundary condition:
\begin{equation}\label{0}
-y''(x)+q(x)y(x)=\lambda y(x), \quad y(0)=y^{(j)}(\pi)=0, \quad j=0,1,
\end{equation}
uniquely determine the potential $q(x),$ see \cite{B}, where also local solvability and, actually, stability
of this inverse problem were established in the class of real-valued $q(x)\in L_2(0,\pi).$ Later on, these
results were refined and generalized to other classes of potentials and boundary conditions \cite{Kar, MO,
FY01, HM04, SavShk, BK19, M, L}. Moreover, there appeared methods, which gave global solution for the inverse
Sturm--Liouville problem as well as for inverse problems for other classes of differential operators (see,
e.g., monographs \cite{FY01, M, L, Yur02 }).

In recent years, there appeared a considerable interest in inverse problems also for nonlocal
Sturm--Liouville-type operators with deviating argument (see \cite{Pik91, AlbHryNizh, FrYur12, VladPik16,
ButYur19, ButPikYur17, Ign18, BondYur18-1, VPV19, DV19, Yur20, SS19, Dur, Yang, Wang, BK, BMSh, Wang21,
ButHu, DB21, DB21-2} and references therein), which are often more adequate for modelling various real-world
processes frequently possessing a nonlocal nature. A special place in this direction is occupied by the
inverse problems for operators with constant delay.

Fix $\nu\in\{0,1\}.$ For $j=0,1,$ consider the boundary value problem ${\cal L}_{\nu,j}(a,q)$ of the form
\begin{equation}\label{1}
-y''(x)+q(x)y(x-a)=\lambda y(x),\quad 0<x<\pi,
\end{equation}
\begin{equation}\label{2}
y^{(\nu)}(0)=y^{(j)}(\pi)=0
\end{equation}
with delay $a\in(0,\pi)$ and a complex-valued potential $q(x)\in L_2(0,\pi)$ such that $q(x)=0$ on $(0,a).$

Denote by $\{\lambda_{n,j}\}_{n\ge1}$ the spectrum of ${\cal L}_{\nu,j}(a,q)$ and consider the following
inverse problem.

\medskip
{\bf Inverse Problem 1.} Given the spectra $\{\lambda_{n,0}\}_{n\ge1}$ and $\{\lambda_{n,1}\}_{n\ge1},$ find
the potential $q(x).$

\medskip
Alternatively, one can consider the cases of Robin boundary conditions
\begin{equation}\label{rob}
y^{(\nu)}(0)-\nu hy(0)=y'(\pi)+H_jy(\pi)=0,\quad j=0,1, \quad h,H_0,H_1\in{\mathbb C}, \quad  H_0\ne H_1,
\end{equation}
which, however, can be easily reduced to conditions (\ref{2}), while all involved coefficients $h$ and $H_0$
along with~$H_1$ are uniquely determined by the two spectra (see \cite{Yur20}).

Various aspects of Inverse Problem~1 were studied in \cite{Pik91, FrYur12, VladPik16, ButYur19, ButPikYur17,
BondYur18-1, %BMSh,
VPV19, DV19, Yur20, Dur, DB21, Ign18, DB21-2} and other works. In particular, it is well known that the two
spectra uniquely determine the potential $q(x)$ as soon as $a\ge\pi/2,$ when the inverse problem is even
overdetermined (see \cite{ButYur19}). For $a<\pi/2,$ the dependence of the characteristic function of any
problem consisting of (\ref{1}) and (\ref{2}) on the potential is nonlinear, which finally results in a
nonlinear equation for recovering $q(x),$ which is referred as the {\it main equation} of the inverse
problem. Specifically, for $a\in[\pi/(N+1),\pi/N)$ with some $N\ge2,$ the main equation will contain
nonlinear integral terms of the form:
\begin{equation}\label{nonlin}
C_k%\underbrace{
\int\limits_{\alpha_{1,k}^0}^{\alpha_{1,k}^1}\ldots\int\limits_{\alpha_{k,k}^0}^{\alpha_{k,k}^1} %}_{k\; {\rm times}}
q(\tau_1)\ldots q(\tau_k)\,d\tau_1\ldots d\tau_k, \quad k=\overline{1,N}, \quad C_N\ne0,
\end{equation}
with certain limits of integration $\alpha_{p,k}^l= \alpha_{p,k}^l(x,\tau_1,\ldots,\tau_{p-1}),$ where $1\le
p\le k\le N$ and $l=0,1.$

Actually, it became a {\it basic question} of the inverse spectral theory for functional-differential
operators with constant delay whether the uniqueness holds also in the nonlinear case~$a\in(0,\pi/2).$ A few
years ago, a positive answer in the case $a\in[2\pi/5,\pi/2)$ was given in \cite{BondYur18-1} for $\nu=0$ and
independently in \cite{VPV19} for $\nu=1.$ Moreover, the overdetermination remains. However, recent authors'
papers \cite{DB21} and \cite{DB21-2} gave a negative answer as soon as $a\in[\pi/3,2\pi/5)$ for the cases
$\nu=0$ and $\nu=1,$ respectively. Specifically, for each such $a$ and $\nu,$ we constructed infinite
families of different iso-bispectral potentials~$q(x),$ i.e. for which both problems consisting of (\ref{1})
and (\ref{2}) possess one and the same pair of spectra. This appeared quite unexpected, in particular,
because of the inconsistence with Borg's classical uniqueness result for~$a=0.$

The general strategy of constructing such iso-bispectral potentials involved establishing a special Fredholm
linear integral operator $M_h$ generated by the restriction $h(x)$ of $q(x)$ to some proper subinterval. Then
an appropriate eigenfunction of $M_h$ was used as a part of the iso-bispectral potentials on some other
proper subinterval, while the variety of the potentials was achieved by multiplying the eigenfunction with
arbitrary complex constants. Meanwhile, the case of a Neumann boundary condition at zero, i.e. when $\nu=1,$
appeared to be more difficult than the case of a Dirichlet one $(\nu=0)$ since the former required finding
such an operator $M_h$ that would possess an eigenfunction with the zero mean value. Even though the
existence of such an operator was not in doubt, finding its concrete example appeared to be a quite difficult
task. After a series of computational experiments we constructed several numerical examples, one of which
fortunately admitted a precise elementary implementation. The numerical simulation allowed us to construct
also elementary iso-bispectral $W_2^1$-potentials for $a\in(\pi/3,2\pi/5)$ and $\nu=0$ (see \cite{DB21-2}).

However, the works \cite{DB21} and \cite{DB21-2} left unclear whether it is also impossible to uniquely
recover the potential by the two spectra in the most difficult nonlinear case $a\in(0,\pi/3),$ when the main
equation would contain nonlinear terms of the form (\ref{nonlin}) with the nonlinearity parameter $N>2,$
becoming unboundedly large for arbitrarily small $a.$ At the same time, the case of small $a$ is especially
interesting and important since it approximates the classical case $a=0,$ in which the uniqueness is well
known.

In the present paper, we address this gap by giving a negative answer also in the remarkable case
$a\in(0,\pi/3)$ for both types of boundary conditions at zero, i.e. for $\nu=0$ and for $\nu=1.$

The paper is organized as follows. In the next section, we study sine- and cosine-type solutions of equation
(\ref{1}) and find some appropriate representations for the characteristic functions of the problems ${\cal
L}_{\nu,j}(a,q)$ for $\nu,j=0,1.$ The main results of the paper are provided in Section~3.
\\

{\large\bf 2. Sine- and cosine-type solutions. The characteristic functions}
\\

For $\nu=0,1,$ denote by $y_\nu(x,\lambda)$ the unique solution of equation (\ref{1}) under the initial
conditions $y_\nu^{(j)}(0,\lambda)=\delta_{\nu,j},$ $j=0,1,$ where $\delta_{\nu,j}$ is the Kronecker delta.
Here and below, $f^{(j)}$ as well as~$f'$ denote the corresponding derivatives of $f$ with respect to the
{\it first} argument. Put $\rho^2=\lambda$ and denote
$$
\omega(x)=\int\limits_a^x q(t)\,dt, \quad \omega_1(x)=\int\limits_{2a}^x q(t)\omega(t-a)\,dt, \quad
y_{0,0}(x,\lambda)=\cos\rho x, \quad y_{1,0}(x,\lambda)=\frac{\sin\rho x}\rho.
$$
The functions $y_0(x,\lambda)$ and $y_1(x,\lambda)$ are called cosine- and sine-type solutions of equation
(\ref{1}), which, in the case of the zero potential, coincide with $y_{0,0}(x,\lambda)$ and
$y_{1,0}(x,\lambda),$ respectively.

For any pair of the parameters $\nu,j\in\{0,1\},$ eigenvalues of the boundary value problem ${\cal
L}_{\nu,j}(a,q)$ with account of their multiplicities coincide with zeros of the entire function
$\Delta_{\nu,j}(\lambda) =y_{1-\nu}^{(j)}(\pi,\lambda),$ which is called {\it characteristic function} of the
problem ${\cal L}_{\nu,j}(a,q).$

For each $\nu=0,1,$ Lagrange's method of variation of parameters leads to the integral equation
$$
y_\nu(x,\lambda)=y_{\nu,0}(x,\lambda) +\int\limits_a^x\frac{\sin\rho(x-t)}\rho q(t)y_\nu(t-a,\lambda)\,dt.
$$
Solving this equation by the method of successive approximations, we arrive at the series
\begin{equation}\label{series}
y_\nu(x,\lambda)=\sum_{k=0}^\infty y_{\nu,k}(x,\lambda), \quad
y_{\nu,k}(x,\lambda)=\int\limits_a^x\frac{\sin\rho(x-t)}\rho q(t)y_{\nu,k-1}(t-a,\lambda)\,dt, \quad k\ge1.
\end{equation}
Clearly, $y_{\nu,1}(x,\lambda)=0$ for $x\in[0,a].$ Moreover, the induction gives $y_{\nu,k}(x,\lambda)=0$ for
$x\in[0,ka]\cap[0,\pi].$ In particular, we have $y_{\nu,k}(x,\lambda)\equiv0$ for $k>N,$ where $N\in{\mathbb
N}$ is such that $a\in[\pi/(N+1),\pi/N).$ Hence, formulae (\ref{series}) take the form
\begin{equation}\label{series-1}
y_\nu(x,\lambda)=\sum_{k=0}^N y_{\nu,k}(x,\lambda), \quad
y_{\nu,k}(x,\lambda)=\int\limits_{ka}^x\frac{\sin\rho(x-t)}\rho q(t)y_{\nu,k-1}(t-a,\lambda)\,dt, \quad
k\ge1.
\end{equation}
The following lemma gives explicit formulae for the terms $y_{\nu,k}(x,\lambda)$ for $\nu=0,1$ and $k=1,2.$

\medskip
{\bf Lemma 1. }{\it For $\nu=0,1,$ the following representations hold:
\begin{equation}\label{y-nu-1}
y_{\nu,1}(x,\lambda)=\frac{\omega(x)}{2(-\lambda)^\nu}y_{1-\nu,0}(x-a,\lambda)
+\frac1{2\lambda^\nu}\int\limits_a^x q(t)y_{1-\nu,0}(x-2t+a,\lambda)\,dt, \quad a\le x\le\pi,
\end{equation}
\begin{equation}\label{y-nu-2}
y_{\nu,2}(x,\lambda)=\frac1{2\lambda^\nu}\int\limits_\frac{3a}2^{x-\frac{a}2}
P_\nu(x,t)y_{1-\nu,0}(x-2t+a,\lambda)\,dt, \quad 2a\le x\le\pi,
\end{equation}
where
\begin{equation}\label{W1}
P_\nu(x,t)=\int\limits_{t+\frac{a}2}^x q(\tau)\,d\tau\int\limits_a^{t-\frac{a}2} q(\xi)\,d\xi
+(-1)^\nu\int\limits_a^{x-t+\frac{a}2} q(\tau)\,d\tau\int\limits_{t+\tau-\frac{a}2}^x q(\xi)\,d\xi, \quad
\frac{3a}2\le t\le x-\frac{a}2\le\pi-\frac{a}2.
\end{equation}}

\medskip
{\it Proof.} Fix $\nu\in\{0,1\}.$ Applying the obvious unified trigonometric relation
\begin{equation}\label{trig}
\frac{\sin\rho(x-t)}\rho y_{\nu,0}(\xi,\lambda)=\frac1{2\lambda^\nu}\Big((-1)^\nu
y_{1-\nu,0}(x-t+\xi,\lambda) +y_{1-\nu,0}(x-t-\xi,\lambda)\Big),
\end{equation}
with $\xi=t-a$ to the second formula in (\ref{series-1}) for $k=1,$ we arrive at representation
(\ref{y-nu-1}).

Further, substituting (\ref{y-nu-1}) into the second formula in (\ref{series-1}) for $k=2,$ we obtain the
relation
$$
y_{\nu,2}(x,\lambda)=\frac1{2\lambda^\nu}\int\limits_{2a}^x\frac{\sin\rho(x-t)}\rho q(t) \Big(
(-1)^\nu\omega(t-a)y_{1-\nu,0}(t-2a,\lambda) +\int\limits_a^{t-a} q(\tau)y_{1-\nu,0}(t-2\tau,\lambda)\,d\tau
\Big)dt.
$$
Using (\ref{trig}) for $\xi=t-2a$ and for $\xi=t-2\tau,$ we get
$$
y_{\nu,2}(x,\lambda)=-\frac{\omega_1(x)}{4\lambda}y_{\nu,0}(x-2a,\lambda) +
\frac{(-1)^\nu}{4\lambda}\int\limits_a^{x-a} q(t+a)y_{\nu,0}(x-2t,\lambda)\,dt \int\limits_a^t q(\tau)\,d\tau
\qquad\qquad\qquad
$$
$$
\qquad-\frac{(-1)^\nu}{4\lambda}\int\limits_{2a}^x q(t)\,dt \int\limits_a^{t-a} q(\tau)
y_{\nu,0}(x-2\tau,\lambda)\,d\tau +\frac1{4\lambda}\int\limits_{2a}^x q(t)\,dt \int\limits_a^{t-a} q(t-\tau)
y_{\nu,0}(x-2\tau,\lambda)\,d\tau.
$$
After changing the order of integrations in the last two terms, we obtain the representation
\begin{equation}\label{y-nu-2-1}
y_{\nu,2}(x,\lambda)=-\frac{\omega_1(x)}{4\lambda}y_{\nu,0}(x-2a,\lambda)
+\frac{(-1)^\nu}{4\lambda}\int\limits_a^{x-a} R_\nu(x,t) y_{\nu,0}(x-2t,\lambda)\,dt,
\end{equation}
where
$$
R_\nu(x,t)=q(t+a)\int\limits_a^tq(\tau)\,d\tau -q(t)\int\limits_{t+a}^xq(\tau)\,d\tau +(-1)^\nu
\int\limits_{t+a}^x q(\tau)q(\tau-t)\,d\tau.
$$
Integrating with respect to the second argument, we get
$$
\int\limits_t^{x-a}R_\nu(x,\tau)\,d\tau=\int\limits_{t+a}^xq(\tau)\,d\tau\int\limits_a^{\tau-a}q(\xi)\,d\xi
-\int\limits_t^{x-a}q(\tau)\,d\tau\int\limits_{\tau+a}^xq(\xi)\,d\xi +(-1)^\nu\int\limits_t^{x-a}d\tau
\int\limits_a^{x-\tau} q(\xi+\tau)q(\xi)\,d\xi.
$$
Changing the order of integrations in the last two terms, we arrive at
$$
\int\limits_t^{x-a}R_\nu(x,\tau)\,d\tau=\int\limits_{t+a}^xq(\tau)\,d\tau\int\limits_a^t q(\xi)\,d\xi
+(-1)^\nu\int\limits_a^{x-t}q(\tau)\,d\tau \int\limits_{t+\tau}^x q(\xi)\,d\xi=P_\nu\Big(x,t+\frac{a}2\Big),
$$
where the function $P_\nu(x,t)$ is determined by (\ref{W1}). In particular, this yields
$$
P_\nu\Big(x,\frac{3a}2\Big)=(-1)^\nu\int\limits_a^{x-a}q(t)\,dt \int\limits_{a+t}^x q(\tau)\,d\tau
=(-1)^\nu\int\limits_{2a}^xq(t)\,dt \int\limits_a^{t-a} q(\tau)\,d\tau =(-1)^\nu \omega_1(x).
$$
Hence, the integration by parts in (\ref{y-nu-2-1}) along with the relations
\begin{equation}\label{der}
y_{\nu,0}'(x,\lambda)=(-\lambda)^{1-\nu} y_{1-\nu,0}(x,\lambda), \quad \nu=0,1,
\end{equation}
gives (\ref{y-nu-2}). $\hfill\Box$

\medskip
{\bf Corollary 1. }{\it The following representations hold:
\begin{equation}\label{y-nu-1'}
y_{\nu,1}'(x,\lambda)=\frac{\omega(x)}2y_{\nu,0}(x-a,\lambda) +\frac{(-1)^\nu}2\int\limits_a^x
q(t)y_{\nu,0}(x-2t+a,\lambda)\,dt, \quad a\le x\le\pi,
\end{equation}
\begin{equation}\label{y-nu-2'}
y_{\nu,2}'(x,\lambda)=\frac{(-1)^\nu}2\int\limits_\frac{3a}2^{x-\frac{a}2}
P_\nu(x,t)y_{\nu,0}(x-2t+a,\lambda)\,dt, \quad 2a\le x\le\pi.
\end{equation}
}

{\it Proof.} By virtue of (\ref{der}) along with the relation $w'(x)=q(x),$ formula (\ref{y-nu-1'}) can be
easily obtained by differentiating (\ref{y-nu-1}). Analogously, differentiating (\ref{y-nu-2}), we get
\begin{equation}\label{y-nu-2''}
y_{\nu,2}'(x,\lambda)=A_\nu(x,\lambda)+B_\nu(x,\lambda) +\frac{(-1)^\nu}2\int\limits_\frac{3a}2^{x-\frac{a}2}
P_\nu(x,t)y_{\nu,0}(x-2t+a,\lambda)\,dt, \quad 2a\le x\le\pi,
\end{equation}
where
$$
A_\nu(x,\lambda)=\frac{(-1)^\nu}2 P_\nu\Big(x,x-\frac{a}2\Big)y_{\nu,0}(2a-x,\lambda), \quad B_\nu(x,\lambda)
=\frac1{2\lambda^\nu}\int\limits_\frac{3a}2^{x-\frac{a}2} \frac\partial{\partial
x}P_\nu(x,t)y_{1-\nu,0}(x-2t+a,\lambda)\,dt.
$$
Since, according to (\ref{W1}), $P_\nu(x,x-a/2)=0,$ we have $A_\nu(x,\lambda)\equiv0$ for $\nu=0,1.$ Further,
we obtain
$$
\frac\partial{\partial x}P_\nu(x,t) =q(x)\left(\int\limits_a^{t-\frac{a}2} q(\tau)\,d\tau
+(-1)^\nu\int\limits_a^{x-t+\frac{a}2} q(\tau)\,d\tau\right).
$$
Thus, making the change of the integration variable $t\to(x+a-t)/2,$ we arrive at
$$
B_\nu(x,\lambda) =\frac{q(x)}{4\lambda^\nu}\int\limits_{2a-x}^{x-2a} G_\nu(x,t,\lambda) \,dt, \quad
G_\nu(x,t,\lambda)=\left(\int\limits_a^\frac{x-t}2 q(\tau)\,d\tau +(-1)^\nu\int\limits_a^\frac{x+t}2
q(\tau)\,d\tau\right)y_{1-\nu,0}(t,\lambda).
$$
It remains to note that for each fixed $x\in(2a,\pi]$ and $\lambda\in{\mathbb C},$ both functions
$G_\nu(x,t,\lambda),$ $\nu=0,1,$ are odd with respect to $t\in(2a-x,x-2a).$ Hence, $B_\nu(x,\lambda)\equiv0$
for $\nu=0,1,$ and we arrive at~(\ref{y-nu-2'}). $\hfill\Box$

\medskip
In our target case $a\in(0,\pi/3),$ besides $y_{\nu,1}(x,\lambda)$ and $y_{\nu,2}(x,\lambda),$ the series in
(\ref{series-1}), generally speaking, possesses also $y_{\nu,k}(x,\lambda)$ for $k=\overline{3,N},$ which
causes the appearance of nonlinear integral terms having the form (\ref{nonlin}) for $k>2$ in representations
of the characteristic functions. This could essentially complicate the study. However, our purposes allow to
zeroize the undesirable terms by imposing additional assumptions on the potential. Namely, after assuming
$q(x)=0$ a.e. on $(3a,\pi),$ we arrive at the representations
\begin{equation}\label{series-2}
y_\nu(x,\lambda)=y_{\nu,0}(x,\lambda) +y_{\nu,1}(x,\lambda) +y_{\nu,2}(x,\lambda), \quad \nu=0,1,
\end{equation}
which allow us to derive convenient representations for the characteristic functions.

\medskip
{\bf Lemma 2. }{\it Let $q(x)=0$ a.e. on $(3a,\pi).$ Then, for $\nu,j=0,1,$ the following representations
hold
\begin{equation}\label{3.1}
\left.\begin{array}{c} %%
\displaystyle \Delta_{\nu,\nu}(\lambda)=(-\lambda)^{\nu}\Big(\frac{\sin\rho\pi}\rho-
\omega\frac{\cos\rho(\pi-a)}{2\lambda}
+\frac{(-1)^\nu}{2\lambda}\int\limits_a^{3a} w_\nu(x)\cos\rho(\pi-2x+a)\,dx\Big),\\[4mm]
\displaystyle \Delta_{\nu,j}(\lambda)=\cos\rho\pi +\omega\frac{\sin\rho(\pi-a)}{2\rho}
+\frac{(-1)^j}{2\rho}\int\limits_a^{3a} w_\nu(x)\sin\rho(\pi-2x+a)\,dx, \quad \nu\ne j, %%
\end{array}\right\}
\end{equation}
where $\omega=\omega(\pi)$ and the functions $w_\nu(x)$ are determined by the formula
\begin{equation}\label{3.5}
w_\nu(x)=\left\{\begin{array}{l}
\displaystyle q(x),\quad x\in\Big(a,\frac{3a}2\Big)\cup\Big(\frac{5a}2,3a\Big),\\[3mm]
\displaystyle q(x)+Q_\nu(x),\quad x\in\Big(\frac{3a}2,\frac{5a}2\Big),
\end{array}\right.
\end{equation}
while
\begin{equation}\label{3.6}
Q_\nu(x)=\int\limits_{x+\frac{a}2}^{3a} q(t)\,dt\int\limits_a^{x-\frac{a}2} q(\tau)\,d\tau
-(-1)^\nu\int\limits_a^{\frac{7a}2-x} q(t)\,dt\int\limits_{x+t-\frac{a}2}^{3a} q(\tau)\,d\tau, \quad
x\in\Big(\frac{3a}2,\frac{5a}2\Big).
\end{equation}
}

\medskip
{\it Proof.} Under the hypothesis of the lemma, formula (\ref{W1}) gives $P_\nu(x,t)=0$ for $t\ge5a/2$ and
$\nu=0,1.$ Thus, substituting $x=\pi$ into (\ref{series-2}) and taking (\ref{y-nu-1})--(\ref{W1}) and
(\ref{y-nu-1'}), (\ref{y-nu-2'}) into account, we arrive at formulae (\ref{3.1}) and (\ref{3.5}) with the
functions $Q_\nu(x),$ $\nu=0,1,$ determined by the formula
$$
Q_\nu(x)=P_{1-\nu}(\pi,x)=\int\limits_{x+\frac{a}2}^\pi q(t)\,dt\int\limits_a^{x-\frac{a}2} q(\tau)\,d\tau
-(-1)^\nu\int\limits_a^{\pi-x+\frac{a}2} q(t)\,dt\int\limits_{x+t-\frac{a}2}^\pi q(\tau)\,d\tau,
$$
where, besides the obvious possibility of replacing both upper limits of integration $\pi$ with $3a,$ it
remains to note that the last internal integral vanishes as soon as $x+t-a/2\ge3a.$ Thus, the upper limit of
integration in the last external integral $\pi-x+a/2$ can be replaced with $7a/2-x,$ which gives (\ref{3.6}).
$\hfill\Box$
\\

{\large\bf 3. The main results}
\\

In this section, we establish non-uniqueness of the solution of Inverse Problem~1 for both cases $\nu=0$ and
$\nu=1$ as soon as $a\in(0,\pi/3).$ Specifically, for each such $\nu$ we construct an infinite family of
iso-bispectral potentials $q(x),$ i.e. for which both problems ${\cal L}_{\nu,0}(a,q)$ and ${\cal
L}_{\nu,1}(a,q)$ possess one and the same pair of spectra. According to the previous section, the spectrum of
any problem ${\cal L}_{\nu,j}(a,q)$ does not depend on $q(x)\in B$ for some subset $B\subset L_2(0,\pi)$ as
soon as neither does the corresponding characteristic function $\Delta_{\nu,j}(\lambda).$ It is sufficient to
restrict oneself with potentials vanishing on the interval $(3a,\pi),$ which allows one to use
representations of the characteristic functions obtained in Lemma~2.

\medskip
{\bf Remark 1.} The difference between the cases $\nu=0$ and $\nu=1$ consists in the following. Since the
functions $\Delta_{\nu,j}(\lambda)$ are entire in $\lambda,$ the first representation in (\ref{3.1}) for
$\nu=0$ implies
\begin{equation}\label{3.6.1}
\omega=\int\limits_a^\pi w_0(x)\,dx,
\end{equation}
which can also be checked directly using (\ref{3.5}) and (\ref{3.6}) for $\nu=0.$ Thus, for $\nu=0,$ the
iso-bispectrality of any subset $B\subset L_2(0,\pi)$ requires only~$w_0(x)$'s independence of $q(x)\in B.$
However, for $\nu=1,$ there is no relation analogous to~(\ref{3.6.1}). In other words, the constant $\omega$
is not determined by $w_1(x).$ Thus, both functions $\Delta_{1,0}(\lambda)$ and $\Delta_{1,1}(\lambda)$ may
depend on $q(x)\in B$ even when $w_1(x)$ does not.

\medskip
Fix $a\in(0,\pi/3)$ and consider the integral operator
\begin{equation}\label{2.1}
M_hf(x)=\int\limits_\frac{3a}2^{\frac{7a}2-x} K_h\Big(x+t-\frac{a}2\Big)f(t)\,dt, \quad \frac{3a}2<x<2a,
\quad {\rm where} \quad K_h(x)=\int\limits_x^{3a} h(\tau)\,d\tau,
\end{equation}
with a nonzero real-valued function $h(x)\in L_2(5a/2,3a).$ Thus, $M_h$ is a nonzero compact Hermitian
operator in $L_2(3a/2,2a)$ and, hence, it has at least one nonzero eigenvalue $\eta.$ Fix $\nu\in\{0,1\}$ and
put $h_\nu(x):=(-1)^\nu h(x)/\eta.$ Then $(-1)^\nu$ is an eigenvalue of $M_{h_\nu}.$ Let $e_\nu(x)$ be a
related eigenfunction,~i.e.
\begin{equation}\label{2.2}
M_{h_\nu}e_\nu(x)=(-1)^\nu e_\nu(x), \quad \frac{3a}2<x<2a.
\end{equation}
Consider the one-parametric family of potentials $B_\nu:=\{q_{\alpha,\nu}(x)\}_{\alpha\in{\mathbb C}}$
determined by the formula
\begin{equation}\label{2.3}
q_{\alpha,\nu}(x)=\left\{\begin{array}{cl}\displaystyle 0,
 &\displaystyle x\in\Big(0,\frac{3a}2\Big),\\[3mm]
\displaystyle \alpha e_\nu(x), &\displaystyle x\in\Big(\frac{3a}2,2a\Big),\\[3mm]
\displaystyle -\alpha K_{h_\nu}\Big(x+\frac{a}2\Big)\int\limits_\frac{3a}2^{x-\frac{a}2} e_\nu(t)\,dt,
 &\displaystyle x\in\Big(2a,\frac{5a}2\Big),\\[3mm]
h_\nu(x), &\displaystyle x\in\Big(\frac{5a}2,3a\Big).
\end{array}\right.
\end{equation}

\medskip
{\bf Lemma 3. }{\it For each fixed $\nu\in\{0,1\},$ the function $w_\nu(x)$ constructed by the formulae
(\ref{3.5}) and (\ref{3.6}) with $q(x)=q_{\alpha,\nu}(x)$ determined by (\ref{2.3}) is independent of
$\alpha.$}

\medskip
{\it Proof.}  Fix $\nu\in\{0,1\}$ and let $q(x)=0$ on $(a,3a/2).$ Then (\ref{3.6}) takes the form
$$%\begin{equation}\label{3.8}
Q_\nu(x)=\left\{\begin{array}{cl} \displaystyle -(-1)^\nu\int\limits_\frac{3a}2^{\frac{7a}2-x}
q(t)\,dt\int\limits_{x+t-\frac{a}2}^{3a} q(\tau)\,d\tau, & \displaystyle x\in\Big(\frac{3a}2,2a\Big),\\[3mm]
\displaystyle \int\limits_{x+\frac{a}2}^{3a} q(t)\,dt\int\limits_\frac{3a}2^{x-\frac{a}2} q(\tau)\,d\tau,
 & \displaystyle x\in\Big(2a,\frac{5a}2\Big),
\end{array}\right.
$$%\end{equation}
which along with (\ref{3.5}) and (\ref{2.1}) implies
\begin{equation}\label{3.9}
w_\nu(x)=\left\{\begin{array}{cl} \displaystyle 0,& \displaystyle x\in\Big(a,\frac{3a}2\Big),\\[3mm]
\displaystyle q(x)-(-1)^\nu M_qq(x), & \displaystyle x\in\Big(\frac{3a}2,2a\Big),\\[3mm]
\displaystyle q(x)+K_q\Big(x+\frac{a}2\Big)\int\limits_\frac{3a}2^{x-\frac{a}2} q(t)\,dt,
 & \displaystyle x\in\Big(2a,\frac{5a}2\Big),\\[3mm]
\displaystyle q(x), & \displaystyle x\in\Big(\frac{5a}2,3a\Big).
\end{array}\right.
\end{equation}
Substituting $q(x)=q_{\alpha,\nu}(x)$ determined by (\ref{2.3}) into (\ref{3.9}) and taking (\ref{2.2}) into
account, we arrive at
$$
w_\nu(x)=\left\{\begin{array}{cl} \displaystyle 0, & \displaystyle x\in\Big(a,\frac{5a}2\Big),\\[3mm]
h_\nu(x), & \displaystyle x\in\Big(\frac{5a}2,\pi\Big),
\end{array}\right.
$$
which finishes the proof. $\hfill\Box$

\medskip
Lemma~3 along with (\ref{3.1}) for $\nu=0$ and relation (\ref{3.6.1}) gives the following theorem.

\medskip
{\bf Theorem 1. }{\it For $j=0,1,$ the spectrum of the problem ${\cal L}_{0,j}(a,q_{\alpha,0})$ is
independent of $\alpha.$ }

\medskip
Thus, $B_0$ is a set of iso-bispectral potentials in the case $\nu=0,$ i.e. Inverse Problem~1 is {\it not
uniquely solvable for $\nu=0$ and $a\in(0,\pi/3).$} However, according to Remark~1, the set $B_1,$ generally
speaking, does not consist of iso-bispectral potentials for $\nu=1.$ Nevertheless, the following theorem
holds.

\medskip
{\bf Theorem 2. }{\it For $j=0,1,$ the spectrum of the problem ${\cal L}_{1,j}(a,q_{\alpha,1})$ is
independent of $\alpha$ as soon~as
\begin{equation}\label{2.3.1}
\int\limits_\frac{3a}2^{2a} e_1(x)\,dx=0.
\end{equation}}

{\it Proof.} According to (\ref{3.1}) for $\nu=1$ and Lemma~3, it is sufficient to prove that the value
$$
\omega=\int\limits_a^\pi q_{\alpha,1}(x)\,dx
$$
is independent of $\alpha$ when (\ref{2.3.1}) holds. Integrating the third line in (\ref{2.3}) for $\nu=1,$
we get
$$
{\cal I}:=\int\limits_{2a}^\frac{5a}2 K_{h_1}\Big(x+\frac{a}2\Big)\,dx\int\limits_\frac{3a}2^{x-\frac{a}2}
e_1(t)\,dt =\int\limits_{2a}^\frac{5a}2 K_{h_1}\Big(x+\frac{a}2\Big)\,dx
\int\limits_\frac{3a}2^{x-\frac{a}2}e_1(x+a-t)\,dt.
$$
Changing the order of integration and then the internal integration variable, we calculate
$$
{\cal I}=\int\limits_\frac{3a}2^{2a}dx \int\limits_{x+\frac{a}2}^\frac{5a}2 K_{h_1}\Big(t+\frac{a}2\Big)
e_1(t+a-x)\,dt = \int\limits_\frac{3a}2^{2a}dx \int\limits_\frac{3a}2^{\frac{7a}2-x}
K_{h_1}\Big(x+t-\frac{a}2\Big) e_1(t)\,dt,
$$
which along with the first equality in (\ref{2.1}) as well as (\ref{2.2}) for $\nu=1$ and (\ref{2.3.1})
implies
$$
{\cal I}= \int\limits_\frac{3a}2^{2a} M_{h_1} e_1(x)\,dx =-\int\limits_\frac{3a}2^{2a} e_1(x)\,dx=0.
$$
Thus, according to (\ref{2.3}) for $\nu=1,$ relation (\ref{2.3.1}) gives
$$
\omega=\int\limits_a^\pi q_{\alpha,1}(x)\,dx=\int\limits_\frac{5a}2^{3a} h_1(x)\,dx,
$$
i.e. the value $\omega$ does not depend on $\alpha,$ which finishes the proof. $\hfill\Box$

\medskip
This theorem reduces the construction of iso-bispectral potentials for the problems ${\cal L}_{1,0}(a,q)$ and
${\cal L}_{1,1}(a,q)$ to the question of finding a function $h(x)\in L_2(5a/2,3a)$ such that the operator
$M_h$ has at least one eigenfunction possessing the zero mean value but related to a nonzero eigenvalue. An
answer to this question is given by the following assertion.

\medskip
{\bf Proposition 1. }{\it Put
\begin{equation}\label{2.3.2}
h_1(x):=\frac{6\pi^2}{a^2} \cos\pi\sqrt{10}\Big(3-\frac{x}a\Big), \quad e_1(x):=\cos\frac{4\pi
x}a-\cos\frac{2\pi x}a.
\end{equation}
Then relation (\ref{2.2}) for $\nu=1$ as well as equality (\ref{2.3.1}) are fulfilled.}

\medskip
{\it Proof.} Let us start with (\ref{2.3.1}), which can be checked by the direct substitution:
$$
\int\limits_\frac{3a}2^{2a} e_1(x)\,dx=\frac{a}{2\pi}\Big(\frac12 \sin\frac{4\pi x}a-\sin\frac{2\pi
x}a\Big)\Big|_{x=\frac{3a}2}^{2a}=0.
$$
Further, according to (\ref{2.1}) and (\ref{2.3.2}), we have
\begin{equation}\label{2.3.3}
M_{h_1}e_1(x)=\frac{3\pi}{a\sqrt{10}}(A_2-A_1),
\end{equation}
where
$$
A_j=2\int\limits_\frac{3a}2^{\frac{7a}2-x}\cos\frac{2\pi jt}a
\cdot\sin\pi\sqrt{10}\Big(\frac72-\frac{x+t}a\Big)\,dt, \quad j=1,2.
$$
For each $j\in\{1,2\},$ we calculate:
$$
A_j=\sum_{k=0}^1\int\limits_\frac{3a}2^{\frac{7a}2-x}\sin\pi\Big(\sqrt{10}\Big(\frac72-\frac{x}a\Big)
-\frac{\sqrt{10}+2j(-1)^k}at\Big)\, dt.
$$
Fulfilling the integration, we arrive at the relation
$$
A_j=\frac1\pi\sum_{k=0}^1 \frac{a}{\sqrt{10}+2j(-1)^k} \cos\pi\Big(\sqrt{10}\Big(\frac72-\frac{x}a\Big)
-\frac{\sqrt{10}+2j(-1)^k}at\Big)\Big|_{t=\frac{3a}2}^{\frac{7a}2-x} \qquad\;\;
$$
$$
\qquad=\frac1\pi\sum_{k=0}^1 \frac{a}{\sqrt{10}+2j(-1)^k}\Big( \cos2\pi j\Big(\frac72-\frac{x}a\Big)
-\cos\pi\Big(\sqrt{10}\Big(2-\frac{x}a\Big)-3j(-1)^k\Big)\Big)
$$
$$
=\frac{a\sqrt{10}}{3\pi}\Big(\cos\pi\sqrt{10}\Big(2-\frac{x}a\Big) -\cos\frac{2\pi jx}a \Big), \quad j=1,2.
\qquad\qquad\qquad\qquad\;\;
$$
Substituting the obtained representations for $A_1$ and $A_2$ into (\ref{2.3.3}), we arrive at (\ref{2.2})
for $\nu=1.$ $\hfill\Box$

\medskip
Theorem~2 and Proposition~1 imply that the family $B_1$ constructed by using the functions $h_1(x)$ and
$e_1(x)$ determined in (\ref{2.3.2}) consists of iso-bispectral potentials for the problems ${\cal
L}_{1,0}(a,q)$ and ${\cal L}_{1,1}(a,q).$ Thus, Inverse Problem~1 is {\it not uniquely solvable also in the
case $\nu=1$ as soon as $a\in(0,\pi/3).$}
\\

{\bf Acknowledgement.} The first author is supported by Project 19.032/961-103/19 of the Republic of Srpska
Ministry for Scientific and Technological Development, Higher Education and Information Society. The second
one is supported by Grants 19-01-00102, 20-31-70005 of Russian Foundation for Basic Research.

\end{document}